\documentclass[12pt]{amsart}

\usepackage{common}
\usepackage{mathrsfs}

\newcommand{\dd}{\m{\bold d}}

\title{unstable classes of metric structures}

\author{Saharon Shelah}
\author{Alexander Usvyatsov}

\address{Saharon Shelah\\
Mathematics Department\\
Hebrew University of Jerusalem\\
91904 Givat Ram, Israel\\
}
\address{ Saharon Shelah \\ Department of Mathematics\\
Hill Center-Busch Campus Rutgers, The State University of New Jersey
110 Frelinghuysen Rd Piscataway, NJ 08854-8019, USA}

\address{Alexander Usvyatsov\\
Mathematics Department\\
Hebrew University of Jerusalem\\
91904 Givat Ram, Israel\\
}

\address{Alexander Usvyatsov\\
    Universidade de Lisboa \\
  Centro de Matem\'{a}tica e Aplica\cc{c}\~{o}es Fundamentais\\
  Av. Prof. Gama Pinto,2\\
  1649-003 Lisboa \\
  Portugal}

\thanks{The authors thank C.\ Ward Henson for stimulating discussions and for suggesting to investigate asymptotic stability. 
We also thank Ilijas Farah for helpful comments and suggestions.}
\thanks{Publication no. 928 on Shelah's list of publications.
    Shelah thanks the Israel Science Foundation
   (Grant 1053/11), and the
   European Research Council (Grant 338821)
   for partial support of this research.
}
\thanks{
    Research partially supported by Marie Sklodowska Curie CIG 321915 "ModStabBan". Usvyatsov thanks the FCT (Funda\cc{c}\~{a}o para a Ci\^{e}ncia e a Tecnologia) for partial support
    of this research:
    FCT grant no. SFRH / BPD / 34893 / 2007; FCT Research Project PTDC/MAT/101740/2008; FCT Research Project PTDC/MAT/122844/2010; Development grant IF/01726/2012. }

\date{\today}

\begin{document}

\begin{abstract}
    We prove a strong non-structure theorem for a class of metric structures with
    an unstable pair of formulae. As a consequence, we show that weak
    categoricity (that is, categoricity up to isomorphisms and not
    isometries) implies several  versions of stability. This is the
    first step in the direction of the investigation of weak
    categoricity and weak stability of metric classes.
\end{abstract}

\maketitle

\section{Introduction and Preliminaries}


\subsection{Introduction}

Categoricity and stability are two if the most basic and important
notions in contemporary model theory. It is therefore not a surprise
that the study of pure model theory of classes of metric structures
has begun with the investigations of these concepts and connections
between them. 

In classical model theory, categoricity investigates the connection between two notions of similarity of structures: isomorphism and elementary equivalence. Recall that two structures are called elementarily equivalent they have the same theory, that is, they satisfy the same (first order) properties. An equivalent definition comes from the Keisler-Shelah Theorem: two structures are elementarily equivalent if and only if they have isomorphic ultra-powers. Clearly two isomorphic structures are elementarily equivalent, but, 
in general, two elementarily equivalent structures $M$ and $N$ are rarely isomorphic. 
One ``silly'' reason for this is that $M$ and $N$ may have different cardinalities. Categoricity in power in a sense means that this is 
the only reason. To be precise, a structure $M$ is called \emph{categorical} in power (cardinality) $\lam$ if there is a \emph{unique} (up to isomorphism) structure $N$ of cardinality $\lam$ which is elementarily equivalent to $M$. In a way, it would be more accurate to say that the (first order) \emph{theory} of $M$ is categorical in $\lam$, or that the \emph{elementary class} generated by $M$ (the class of all structures elementarily equivalent to $M$) is.

The study of categoricity led to many important discoveries in model theory and generated a large body of research. The first important milestone was Morley's Theorem \cite{Mor} (affirming \L o\'{s}'s Conjecture) stating that a countable first order theory, categorical in some uncountable cardinal is categorical in every such cardinal. One important step in Morley's proof was showing that uncountable categoricity implies $\aleph_0$-stability (which means that there are only countably many types over any countable set), 
and then exploring properties of $\aleph_0$-stable theories (such as existence of prime models over any set). The main tool in the proof of $\aleph_0$-stability was the \emph{Ehrenfeucht-Mostowski construction} which yields a model $M\models T$ ``quite close'' to a particular ``skeleton'' (in the most basic case, a well-ordered indiscernible sequence $I$) for an arbitrary theory $T$. Such a model is normally called the Ehrenfeucht-Mostowski model (or just the EM-model) for $T$ over $I$. Morley observed that such EM-models always realize only countably many types over all countable sets, concluding that if $T$ is not 
$\aleph_0$-stable, it has at least two non-isomorphic models of any uncountable cardinality (one that realizes uncountably many types over some countable set, and one which does not). One can perhaps regard this result as the first general \emph{non-structure} theorem in model theory. 

The first author of this article generalized this theme in several different directions. The main generalization is concerned with replacing \L o\'s's question ``When is a theory $T$ categorical in big enough cardinalities?'' with a much more general question ``When does a theory $T$ have 'few' models in big enough cardinalities?'', or at least ``When can large enough models of $T$ be described by a small collection of invariants?''. This led to the development of a very broad line of research called ``Classification Theory'' \cite{Sh:c}. The premise of classification theory is that there exist true dichotomies, called ``dividing lines'' so that models of theories that lie on the ``good'' side of a dividing line can be analyzed in some way (and, more generally, possess many ``nice'' properties), whereas the theories on the ``dark side'', that is, on the other side of the dividing line, admit various kinds of non-structure. In particular, going back to the question of counting models up to isomorphism, classification theory promotes the search of properties so that if $T$ has such a property, one can prove various structural results about it, and, ideally, classify its models in terms of a well-understood set of invariants; and if $T$ does not have it, then one can prove a non-structure theorem about it - in this case, that $T$ has ``many'' - ideally, the largest possible number - of non-isomorphic models in some (ideally, every) big enough cardinality.

Stability, which is a generalization of Morley's $\aleph_0$-stability, is one such dividing line, which, perhaps, still remains the major one and the most extensively studied. One way to define stability is the following: a theory $T$ is stable if it is $\lam$-stable for some cardinality $\lam$ (hence generalizing the definition of $\aleph_0$-stability directly). However, an equivalent definition (which turns out to be more convenient for our purposes) is that $T$ is stable if no formula exhibits the \emph{order property}. We shall discuss different analogues of the order property (appropriate for our setting) at length later. 

The ``positive'' side of stability is not enough for classifying models up to isomorphism (this is why stronger properties, such as super-stability, NDOP, NOTOP, etc, were introduced). However, models of stable theories admit many good structural properties, investigating which led to many important discoveries and applications. On the other hand, the non-structure side is very clear: an unstable theory $T$ (say, a countable one) has the maximal possible number of non-isomorphic models, even slightly saturated ones, in any uncountable cardinality. The proof involves the use of ``generalized EM-models'', which is yet another generalization of Morley's techniques. This general technique of proving non-structure has paved the way to a broader subject of ``non-structure tehory'' \cite{Sh:e}. We shall use ideas and basic results from \cite{Sh:e} extensively in this paper. 

Another way in which Morley's work was generalized by the first author, was taking classification theory, and in particular, the study of categoricity, beyond the scope of elementary classes. Techniques used in this paper are general enough to study, for example, classes of elementary substructures of a big homogeneous model $\fC$. This way one expands the discussion to classes of models that omit a given (``nice enough'') collection of types over the empty set.

\smallskip

Our current goal is to generalize results of the kind mentioned above to the setting of metric structures. In particular, we set up a framework for proving non-structure theorems in this context. 

\smallskip

While the notion of categoricity have a very precise meaning in classical model theory, there are
more than one natural choices when topology is introduced as an
inherent part of the structure. For example, while classical
categoricity means having a unique model up to isomorphism, when
working with metric structures, one has to specify which
isomorphisms are being considered.

Uncountable (more precisely, non-separable) 
categoricity of metric classes with respect to \emph{isometries} has
been studied independently in different frameworks in
\cite{Ben05} and \cite{ShUs837}. For a model theorist this is
perhaps the most natural choice of isomorphism, since it preserves
\emph{the whole structure}.  The notion of isomorphism 
considered in 
\cite{ShUs837} is as a matter of fact a bit more general, and more recent works of Hirvonen and Hyttinen \cite{HH} weakened 
the definition further, allowing ``almost isometries'', that is, arbitrarily small perturbations. 

However, for many 
natural classes, even this choice would be considered too strict (too rigid). For example, the
``correct'' notion of isomorphism between Banach spaces is widely
accepted to be a linear homeomorphism, or, equivalently, a $c$-isomorphism: a linear isomorphism $T$
satisfying $\frac 1c \le \|T\| \le c$ for some $c \ge 1$. Dealing with such maps that are allowed
to perturb the structure is significantly more challenging
from the model theoretic point of view. 

Naturally, searching for dividing lines which would lead to classification of metric structures up to topological isomorphism, 
one needs to introduce new properties that take topology into account. It is our thesis that stability will still play an important role in this line of research. However, the notion of stability has to be weakened and made more flexible. In this article we will consider several such possible weakenings. It is not completely clear which one is the ``correct'' one. Most likely, each one of the properties discussed below, is the right notion for a certain structural context. We are not concerned 
with structure results in this article. Our goal here is to discuss the ``negative side'' of the proposed dividing lines.



\smallskip


More precisely, all results in this paper can be summarized as  ``various
notions of instability lead to non-structure'', i.e. to the
existence of many non-isomorphic models. This line of research serves several purposes. 


The first one should already be clear from the discussion above. Suppose we have a metric structure $X$ which does not admit this kind of non-structure. One particular 
strong (and interesting) assumption in this direction is \emph{non-separable categoricity} with respect to 
isomorphisms. Then by the results in this paper, $X$
admits several (weak, topological) versions of stability. This clearly constitutes a first step in the direction of understanding non-separable categoricity in this category. 

Another motivation for our investigation is the search for natural weakening of stability in the metric context. The premise 
of this article is that even when a theory
(or a formula/predicate) lacks the full strength of stability, there
are natural weaker properties (which would not make sense in
classical model theory since their definitions rely heavily on the
existence of a metric) that might hold. We introduce and discuss
several such notions. Some of them have already been investigated in
the context of Banach spaces. For instance,
``Junge-Rosenthal''-stability, which we define in section 4, is
motivated by the notion of an ``asymptotically symmetric'' Banach
space, originally due to Junge and Rosenthal, studied in
\cite{JKO}.

Stability is an important model theoretic property that has already been shown to have various 
applications in many different contexts in mathematics, including functional analysis and Banach space geometry. For example. investigation of quantifier-free
stability for Banach spaces led to the spectacular result of Krivine
and Maurey \cite{KM} on the existence of an almost isometric copy of
$\ell_p$ in any stable Banach space, a line of research pursued further by Iovino (e.g. \cite{Iov:def}). Stability was 
generalized to arbitrary definable predicates by Ben Yaacov and the second author \cite{BU}, and used by the authors of this 
article to establish Henson's Conjecture in \cite{ShUs}. It is therefore natural to conjecture that at least some of the weaker notions of ``metric'' stability introduced in this paper will also have useful properties, both for model theoretic study of classes of classes of metric structures, and for the study of the local and global geometry of such structures.  In this way, this paper opens many doors for future research. 

Besides this paper, there have been not many steps have been made in the direction of non-structure  in the general context of metric structures). One related work is Farah and Shelah \cite{FaSh}, where a strong dichotomy result for a number of ultrapowers of a metric structure of cardinality at most the continuum is shown. In particular, it is shown that instability implies that this number is $2^{2^\aleph_0}$.  On the one hand, here we prove non-structure results in a much more general context: we do not assume that the class is elementary or that the vocabulary is countable; in addition, our results hold in all cardinalities above the cardinality of the language (and not just the continuum). Models that we obtain are different from each other in a strong way, not just non-isometric, and in each context we  explicitly compute a constant $c$, up to which none of the two of these models can be isomorphic  (and moreover, none of them can be $c$-embedded into any other). On the other hand, the models constructed in \cite{FaSh}  are all ultrapowers of one particular structure, which is a stronger non-structure result.

As one example of a consequence
of our analysis, let us include here the following result in model
theory of Banach spaces, which (we believe) should have interesting
applications, such as throwing light on model theoretic properties
of particular Banach spaces (it is proven in section 4):

\begin{thm}\label{thm:assym}
    Let $\mathcal B$ be a Banach space whose continuous theory is
    categorical (with respect to isomorphisms of Banach spaces) in
    some uncountable density. Then it is asymptotically symmetric
    (as defined in \cite{JKO}).
\end{thm}

\subsection{Preliminaries}


The proofs in this article rely mostly on Ehrenfeucht-Mostowski
constructions, hence do not require any compactness assumptions and
can be carried out in very general settings, such as abstract metric
classes (as defined in \cite{ShUs837}) with arbitrary
large models. For the sake of simplicity of presentation, we have decided to focus on the well-developed context of 
a monster metric space (it is the central concept in \cite{ShUs837}), which is
already quite general. 


Since the definition of a monster metric space is rather long, we
will not include it here (the reader is welcome to consult \cite{ShUs837}). However, to make the presentation more self-contained, we will define a simplified version, which we call a homogeneous metric
monster.




Fix a ``big enough'' cardinal $\lam^* = {(\lam^*)}^{\aleph_0}$. We refer to \cite{ShUs837} for the definition of a \emph{metric structure}.

\begin{dfn}\label{dfn:momspace}
  A \emph{homogeneous metric monster} is a metric structure
 of cardinality $\lam^*$ which is strongly homogeneous
  in the following sense: every partial function $f \colon A \to \FC$ (with
  $A \subseteq \FC$ of cardinality less than $\lam^*$) which preserves positive
  existential formulae (see \cite{ShUs837}, Definition 2.13), can be extended to an automorphism of $\fC$.
\end{dfn}

\begin{rmk}
  If $\FC$ is a homogeneous metric monster then it is a momspace (monster
  metric space) as defined in \cite{ShUs837}, Definition
  2.17 for $\Delta(\fC) = $positive existential formulas.
\end{rmk}

From now on, $\FC$ will denote a monster metric space (unless stated
otherwise). The reader can assume without almost any loss of
generality that it is a homogeneous metric monster. All elements,
tuples, sets, will be taken inside $\fC$. The finite diagram of
$\FC$ will be denoted by $D$. By ``types'' we will always mean
$D$-types,
which in case of a homogeneous monster is just the collection of all
positive existential types. By ``formula'', we always mean a
$\Delta(\fC)$-formula, which in case of a homogeneous monster is
simply a positive existential formula.

\medskip

Let $K$ be the class of all almost elementary submodels of the monster $\fC$, $K^c$ the class
of complete such models (that is, almost elementary submodels of $\fC$ which are complete with respect to the metric of $\fC$; see also  \cite{ShUs837}, Definition 2.19). We will call elements of
$K$ \emph{premodels} and elements of $K^c$ \emph{models}.


We recall the definition of approximations of formulae and types
(see more in \cite{ShUs837}, section 3):

\begin{dfn}\label{MOM.2}
\begin{enumerate}
\item
For a formula $\varphi(\bar
x)$ possibly with parameters, we define
$\varphi^{[\varepsilon]}(\bar x) = \exists \bar x'(\varphi(\bar x') \wedge
\bold d(\bar x,\bar x') \le \varepsilon)$.  So
$\varphi^{[\varepsilon]}(\bar x,\bar a) = (\exists \bar x')[\bold
d(\bar x,\bar x') \le \varepsilon \wedge \varphi(\bar x',\bar a)]$.
\item
For a type $p$, define $p^{[\varepsilon]} =
\{(\bigwedge_{\ell < n} \varphi_\ell)^{[\varepsilon]}:n < \omega$ and
$\varphi_0,\ldots,\varphi_{n-1} \in p\}$.
\end{enumerate}
\end{dfn}

Recall that a pair $(\ph(\x),\ps(\x))$ is called
$\eps-contradictory$ (where $\eps>0$ as usual) if
$\dd(\ph^\fC,\ps^\fC) > \eps$. A pair of formulae (maybe with
parameters) $(\ph(\x),\ps(\x))$ is called \emph{contradictory} if it
is $\eps$-contradictory for some $\eps>0$.

\medskip

Since the most important context (or at least the context of interest to most readers) in which our results apply is the context 
of a fixed continuous first order theory $T$ \cite{BU,BBHU},  let us explicitly note the following obvious (but important) fact:

\begin{rmk}
	If $T$ is a complete continuous first order theory, then, naturally, the monster model of $T$ (a ``big enough'' 
	saturated model of $T$) is a homogeneous metric monster. In this case, pre-models of $T$ are precisely the elements of $K$, and 
	(complete) models of $T$ are precisely the elements of $K^c$. It is an easy exercise to relate local (in)stability as defined in \cite{BU} to various versions of instability discussed in this article. 
\end{rmk}








\section{Ehrenfeucht-Mostowski models}


Following comments and questions from several readers who are not model theorists, and in order to make the paper more self-contained, in this section we explain in some detail the construction of Ehrenfeucht-Mostowski models, in particular, the procedure of expanding the language with Skolem functions. Everything here is standard and well known, but we thought that it might be helpful to collect all these facts  in one place, and to explain how exactly they apply in our context. A reader who is familiar with section 5 of \cite{ShUs837}, can safely skip the rest of this section. 


\subsection{Ehrenfeucht-Mostowski construction}



As usual, $\fC$ denotes a homogeneous metric monster in vocabulary $\tau$.

First we expand $\fC$ with Skolem functions. More specifically, for every formula $\ph(x,\y)$, we add to $\tau = \tau(\fC)$ a function symbol $f(\y) = f_{\ph(x,\y)}$. Call the new vocabulary $\hat{\tau}$; so
\[
	{\tau'} = \tau\cup \set{f(\y)_{\ph(x,\y)} \; \colon \;\ph(x,\y) \text{ is a } \tau-\text{formula}}
\]
Note that in the formula $\ph(x,\y)$ above the variable $x$ is a singleton (whereas $\y$ can be a tuple).

We then expand $\fC$ to a $\tau'$-structure ${\fC'}$ so that the following condition holds: for every $\tau$-formula $\ph(x,\y)$ and every $\b \in \fC$, letting $f = f_{\ph(x,\y)}$, if $\fC \models \exists x \ph(x,\b)$, then in particular $\fC\models \ph(a,\b)$ where  $a = f^{{\fC'}}(\b)$. 

\smallskip

Let $A$ be any subset of $\fC$ (hence of ${\fC'}$), and let $M'$ be the $\tau'$-substructure generated by $ A$ in $\\fC'$. Let $M$ be the $\tau$-structure which is the restriction of $M'$ to $\tau$. Then clearly $M$ is a substructure of $\fC$; but since $\tau'$ has Skolem functions for $\tau$, $M$ is in fact an elementary substructure of $\fC$. Hence the metric closure of $M$, $\overline{M}$, is an almost elementary substructure of $\fC$. In particular, $M \in K$ and $\overline{M} \in K^c$; that is, $M$ is a pre-model, and $\overline{M}$ is a model.

Note that ultimately we are interested only in $\tau$-structures, which are restrictions of $\tau'$-structures. Being such a restriction ensures that $M$ is an elementary submodel of $\fC$ (hence $\overline{M}$ is a model). However, we do not care about any other properties of the structure $M'$. In particular, it does not need to be a metric structure. Therefore, we do not impose any conditions on the connections between the expanded language and the metric (e.g., continuity of the Skolem functions). 

\medskip

Let $M$ be a $\tau$-structure, $I$ be an arbitrary infinite linear order, and $\a_I = \inseq{\a}{i}{I}$ an indiscernible sequence in $M$. Then the \emph{EM-type} of $\inseq{\a}{i}{I}$ is the following collection of $n-types$ over $\emptyset$ (in $M$):

\[
	\Phi = \left\{ \tp(\a_{i_0},\ldots, \a_{i_{n-1}}) \colon n<\om, \; i_0<i_1<\ldots<i_{n-1} \right\}
\]

Such a collection $\Phi$ is called an \emph{Ehrenfeucht-Mostowski blueprint}, or just an \emph{EM-blueprint}. Note that by indiscernibility, $\Phi$ contains exactly one $n$-type for each $n<\om$. 

Motivated by this example, we now define EM-blueprints more generally. 

\begin{dfn}\label{dfn:EM-blueprint}
  \begin{enumerate}
  \item
  
    An \emph{Ehrenfeucht-Mostowski blueprint} in vocabulary $\tau$ (a $\tau$-EM bluebrint) is a collections of types $\Phi$ satisfying:
    \begin{enumerate}
      \item $\Phi = \set{p_n\colon n<\om}$
      \item $p_n = p_n(\x_0,\ldots,\x_{n-1})$ is an $n$ complete $\tau$-type (over $\emptyset$).
      \item For every $m<n$ and every $0\le i_0<i_1<\ldots<i_{m-1}\le n$, the restriction of $p_n$ to the variables $i_0, \ldots, i_{m-1}$
      is $p_m$. 
    \end{enumerate}
  \item
  	 Given two vocabularies $\tau \subseteq \tau'$, we say that a $\tau'$-EM-blueprint $\Phi' = \set{p'_n\colon n<\om}$ \emph{implies} or \emph{extends} the $\tau$-blueprint $\Phi = \set{p_n\colon n<\om}$ if $p_n \subseteq p'_n$ for all $n$.
  
  \end{enumerate}
\end{dfn}

\medskip

Coming back to our setting, recall that $\fC'$ is an expansion of $\fC$ with Skolem functions. Let $\Phi$ be a $\tau'$-EM-blueprint, $J$ an order type such that there exists an indiscernible sequence $\a_J = \inseq{\a}{i}{J}$ in $\fC'$ whose EM-type is $\Phi$. Then we denote by $EM(J,\Phi)$ the functional closure of $\a_J$ in $\fC'$, and let $EM_\tau(J,\Phi)$ be the $\tau$-reduct of $EM(J,\Phi)$. Since $\fC'$ has Skolem functions, clealry (as above) $EM_\tau(J,\Phi) \prec \fC$. In other words,  $EM_\tau(J,\Phi) \in K$, and $\overline{EM_\tau(J,\Phi)} \in K^c$. This gives us a method of constructing a model (in the metric sense) from an indiscernible sequence. 

Of course there may be many sequences of order type $J$ with EM-type $\Phi$ in $\fC'$, hence technically, the EM-structures discussed above are determined only up to isomorphism. In other words, $EM(J,\Phi)$ really refers to the isomorphism type of the functional closure of any sequence of EM-type $\Phi$ in $\fC'$; similarly for $EM_\tau(J,\Phi)$. We will nevertheless often talk about these as if they are specific structures, and write things like ``$EM_\tau(J,\Phi) \in K$''; but no confusion should arise. 

%
%

\begin{dfn}
\begin{enumerate}
\item
We call $\Phi$ \emph{proper} for $K$ (or for $\fC$) if for any order type $J$, the structure $\EM_\tau(J,\Phi)$ is in $K$.
\item
We call a $\tau$-EM-blueprint $\Phi_0$ \emph{proper} for $K$ (or for $\fC$) if there exists a $\tau'$-blueprint $\Phi$ which extends $\Phi_0$ and is proper for $\fC$.
\end{enumerate}
\end{dfn}


It follows from the discussion above that $\Phi$ is proper for $\fC$ if for every order type $J$, $\fC'$ contains an indiscernible sequence of order type $J$ whose EM-type is $\Phi$.

\bigskip

\subsection{Existence of proper EM-blueprints} 
\label{sub:existence_of_proper_em_blueprints}

Let $T$' be a discrete (classical) first order theory in a vocabulary $\rho'$, $\Gamma$ a collection of types in $S(T')$. Recall that $EC(T',\Gamma)$ denotes the collection of all models of $T'$ omitting all types in $\Gamma$.  Given a vocabulary $\rho \subseteq \rho'$, we denote by $PC_\rho(T',\Gamma)$ the collection of all $\rho$-reducts of models in $EC(T',\Gamma)$ ($PC$ stands for ``Pseudoelementary Class'').

Let us recall the following strengthening (due to the first author) of Morley's Omitting Types Theorem \cite{MorOT},\cite{Sh:c}(Ch.VII). We refer the reader to Appendix A of Baldwin \cite{Bal-book} for a concise presentation, definitions, and proof. We only state here the parts that are relevant for our purposes. 

\begin{thm}\label{thm:MorleyOT}(Generalized Morley's Omitting Types Theorem)
	Let $T$ be a discrete $\rho$-theory with Skolem functions, $\Gamma$ a collection of partial $\rho$-types (in finitely many variables)  over $\emptyset$, and let $\mu = (2^{\tau})^+$. 
	\begin{enumerate}
	\item
	Suppose that for every $\al<\mu$ 	there is $M_\al \in EC(T,\Gamma)$ such that $|M_\al|>\beth_\al$. 
	Then there exists an EM-blueprint $\Phi$ realized in every $M_\al$ proper for $EC(T,\Gamma)$. 
	
	In other words, there exists an EM-blueprint $\Phi$ such that:
	\begin{itemize}
		\item In every $M_\al$ there exists an indiscernible sequence \lseq{\a}{i}{\om} whose EM-type is $\Phi$.
		\item For every order type $J$ we have $EM(J,\Phi) \in EC(T,\Gamma)$.
	\end{itemize}
	
	\item Suppose $P$ is a one-place $\tau$-predicate. If, in addition to the assumption of (i), we have $|P^{M_\al}|\ge \beth_\al$, then 
	the above indiscernible sequence can be chosen ``in $P$''. 
	
	In other words, we can find $\Phi$ as in clause (i) such that the formula $P(\x_0)$ is in $p_1(\x_0)$ (and 
	hence $P(x_i)$ is implied by $p_n(x_0, \ldots, x_{n-1}) $ for all $n$ and $i<n$). 
	\end{enumerate}

\end{thm}

\smallskip

Let us go back to our context, that is,  $\fC$ is a homogeneous metric monster, $\fC'$ is an expansion of $\fC$ to a vocabulary $\tau'$ with Skolem functions. 

\begin{cor}\label{cor:EM-exist}(Existence of proper EM blueprints).
Let $\Phi$ be an EM-blueprint in vocabulary $\tau$. Let $\mu = (2^{\tau})^+$ be as in Theorem \ref{thm:MorleyOT}. Assume now that for any $\al<\mu$ there is an indiscernible sequence $\a_{I_\al}$ of order type $I_\al$, $|I_\al|\ge\beth_\al$ in $\fC$, whose EM blueprint is $\Phi$. Then $\Phi$ is proper for $\fC$.
\end{cor}
\begin{prf}
	Expand $\tau'$ with a new one-place predicate symbol $P$, and further expand the new vocabulary with Skolem functions; call the new vocabulary $\hat \tau$. For every $\al <\mu$, let $M_\al$ be the functional closure of $I_\al$, so that $P$ is interpreted as $I_\al$. 
	
	Since $\fC$ is a homogeneous monster, there exist (see e.g. \cite{Sh:3}) a discrete complete first order $\tau$-theory $T$, and a collection $\Gamma$ of $\tau$-types over the empty set, such that $K = EC(T,\Gamma)$. Of course, $\Gamma$ is a collection of \emph{partial} types in $\hat\tau$, and $T$ is an $\hat\tau$-theory. Clearly, each $\hat M_\al$ is therefore a model in $EC(T,\Gamma)$, and together they satisfy the assumptions of Theorem \ref{thm:MorleyOT}. Hence there exists an $\hat\tau$-EM-blueprint $\hat \Phi$ as there. In other words, $EM(J,\hat \Phi)$ is a $\hat \tau$-structure, which is a model of $T$ and omits all types in $\Gamma$; hence $EM_\tau(J,\hat \Phi) \in EC(T,\Gamma)$ (for any order type $J$). Moreover, by clause (ii) if the Theorem, $\hat \Phi$ can be chosen ``in $P$''. In particular, $\hat \Phi$ extends $\Phi$.

	Restricting ourselves to the vocabulary $\tau'$, and the EM-blueprint $\hat\Phi$ to $\Phi'$, we obtain a $\tau'$-EM-blueprint, which is clearly still proper for $EC(T,\Gamma) = K$, and it still extends $\Phi$, which is exactly what we wanted. 
\end{prf}

\begin{rmk}
	A slight modification of the proof above works for any pseudoelementary class $PC_\tau(T,\Gamma)$, hence in essence for any abstract elementary class $K$. Since all that is used in the proofs in this paper is instability and the existence of proper EM-blueprints, analogous results can be obtained -- by essentially the same arguments -- in much more general contexts (e.g., abstract metric classes with arbitrarily large models). However, we will not go into details of this here. The main motivation for our work is continuous elementary classes, and we feel that by discussing homogeneous monsters, we have already gone far enough outside of this ``comfort zone''. Nevertheless, we preferred to exemplify the general method, rather than rely on compactness, which is unnecessary for the type of non-structure theorems that are shown here. 
\end{rmk}



\section{Instability - the  general case} 
\label{sec:instability}

As before, $\fC$ is a homogeneous metric monster. 

\begin{dfn}\label{dfn:unstable}
\begin{enumerate}
\item
    A contradictory pair $(\ph(\x,\y),\ps(\x,\y))$ is called
    \emph{unstable} if there are arbitrary long indiscernible
    sequences $\lseq{\a}{i}{\lam}$, $\lseq{\b}{i}{\lam}$ such that
    $$i<j<\lam \then \ph(\a_i,\b_j), \ps(\a_j,\b_i)$$
    In this case we also say that the pair above has the \emph{order
    property}.
\item
    A pair as above is called $c$-additively unstable
    or $c$-add-unstable (for $c \ge 0$) if in addition
    $\dd(\ph,\ps) \ge c$.
\item
    $\fC$ is called unstable if there exists an
    unstable pair.
\item
    $\fC$ is called $c$-additively unstable
    or $c$-add-unstable if there exists a
    $c$-add-unstable pair.
\end{enumerate}
\end{dfn}

Clearly, if $\fC$ is compact (see \cite{ShUs837}), e.g., is the monster model of a continuous first order theory, then the pair $(\ph(\x,\y),\ps(\x,\y))$ is unstable if and only if there exists an infinite indiscernible sequence (of any order type) as in Definition \ref{dfn:unstable}(i). But since we have discussed the general theory of proper EM-blueprints in the previous section, we can now make a more interesting observation:

\begin{obs}\label{obs:unstable}
	\begin{enumerate}
	\item
	Let $\fC$ be a homogeneous metric monster in vocabulary $\tau$, $\mu = (2^{|\tau|})^+$, and let $(\ph(\x,\y),\ps(\x,\y))$ be a contradictory pair. Assume that for each $\al<\mu$ there are indiscernible sequences $\inseq{\a}{i}{I_\al}$, $\lseq{\b}{i}{I_\al}$ witnessing instability of $(\ph(\x,\y),\ps(\x,\y))$ (as in Definition \ref{dfn:unstable}(i)) with $|I_\al|\ge\beth_\al$.
	
	Then the pair $(\ph(\x,\y),\ps(\x,\y))$ is unstable. 
	\item
		Moreover, under the assumptions of clause (i), there exists an EM-blueprint $\Phi$ proper for $\fC$ such that
    if $I$ is a linear order and $M = \EM(I,\Phi)$ with skeleton
    \inseq{\a_i\b}{i}{I}, then in $M$ we have
\begin{itemize}
\item
    $i<j \in I \then \ph(\a_i,\b_j)\land\psi(\a_j,\b_i)$
\end{itemize}

	\item
		Analogous results hold for $c$-additively unstable.
	\end{enumerate}
\end{obs} 
\begin{prf}
	Clause (i) follows from clause (ii), which, in turn, follows from Corollary \ref{cor:EM-exist}. Clause (iii) is similar. 
\end{prf}


\section{Instability implies strong nonstructure}






The following is a simplification of Definition III.2.5(1) in
\cite{Sh:e}:

\begin{dfn}
\begin{enumerate}
\item
    Let $I$ be a linear order, $\Phi$ an EM-blueprint proper for $\fC$, $M = \EM(I,\Phi)$,
    $\inseq{\a}{i}{I}$ is the skeleton of $M$, $\c, \d \in M$. We
    say that $\c$ and $\d$ are \emph{similar}, $\c\sim\d$ if
    \begin{itemize}
    \item
        $\c = \bar\sigma(\a_{\bar i}), \d = \bar\sigma(\a_{\bar j})$
    \item
        $\bar i$ and $\bar j$ have the same quantifier free type in
        $I$ in the language of order \set{<}.
    \end{itemize}
\item
    Let $\theta(\x,\y)$ be a formula in the language or order $\set{<}$, and let $I$, $J$
    be linear orders. We say that $I$ is $\theta(\x,\y)$-unembeddable
    into $J$ if for any EM-blueprint $\Phi$ proper for $\fC$ and a function $f \colon I \to \EM(J,\Phi)$
    there exist $\bar i_1,\bar i_2 \in I$ such that $I \models
    \theta(\bar i_1,\bar i_2)$ and $f(\bar i_1) \sim f(\bar i_2)$.
\end{enumerate}
\end{dfn}

\begin{fct}\label{fct:unembed}
    Let $\theta(x_0x_1,y_0y_1)$ be the following formula in the language of
    order: $\theta = x_0<x_1\land y_1<y_0$. Let
    $\lam>|L|$. Then there exists a family $\lseq{I}{\al}{2^\lam}$
    of linear orders of cardinality $\lam$ which are pairwise
    $\theta(\x,\y)$-unembeddable into each other.
\end{fct}
\begin{prf}
    Combine
    \cite{Sh:e} VI.3.1(2) and III.2.21
    (see also Definition
    III.2.5(2) there).
\end{prf}

%
We now define weak embeddings between metric structures.

\begin{dfn}
\begin{enumerate}
\item
    Let $\Delta$ be a set of formulae, $c \ge 0$, and let $M,N
    \in K$. We say that $f \colon M \to N$ is a
    \emph{$(\Delta,c)$-additive-embedding}
    or $(\Delta,c)$-add-embedding if it is 1-1 and for every $\ph(\x) \in
    \Delta$, $\a \in M$ we have $M \models \ph(\a) \then N \models
    \ph^{[c]}(f(\a))$.
\item
    We say that $f\colon M\to N$ is a \emph{$(\Delta,c)$-additive
    isomorphism} or $(\Delta,c)$-add-isomorphism if it is an onto
    $(\Delta,c)$-embedding.
\end{enumerate}
\end{dfn}

\begin{prp}\label{prp:unstable}
    Denote $\theta(x_0x_1,y_0y_1) = x_0<x_1 \land y_1<y_0$ and
    let $I$, $J$ be linear orders such that $I$ is
    $\theta$-unembeddable into $J$.
    Assume that $\fC$ is \emph{unstable} with respect to some
    $c$-contradictory pair $(\ph(\x,\y),\psi(\x,\y))$.
    Then there exists an EM-blueprint $\Phi$ such that there is no
    $(\set{\ph(\x,\y),\psi(\x,\y)},\frac c 3)$-add-embedding $f \colon \overline{\EM(I,\Phi)} \to
    \overline{\EM(J,\Phi)}$.
\end{prp}
\begin{prf}
    By instability of the pair $(\ph,\psi)$ and the standard
    Erd\"os-Rado argument


\begin{tens_equ}\label{equ:Phi}
    There exists an EM-blueprint $\Phi$ proper for $\fC$ such that
    if $I$ is a linear order and $M = \EM(I,\Phi)$ with skeleton
    \inseq{\a_i\b}{i}{I}, then in $M$ we have
\begin{itemize}
\item
    $i<j \in I \then \ph(\a_i,\b_j)\land\psi(\a_j,\b_i)$
\end{itemize}
\end{tens_equ}


    Denote $M = \overline{\EM(I,\Phi)}$, $N =
    \overline{\EM(J,\Phi)}$ and let $f \colon M \to N$.
    Assume by contradiction that $f$ is a $(\set{\ph,\psi},\frac c 3)$-embedding.

    Let $\inseq{\a_i\b}{i}{I}$, $\inseq{\a'_j\b'}{j}{J}$ be the
    skeletons of $\EM(I,\Phi),\EM(J,\Phi)$ respectively.
    Since $I$ is $\theta$-unembeddable into $J$ (and we can identify
    $I$ with $\inseq{\a_i\b}{i}{I}$), there exist
    $i_0,i_1,l_0,l_1 \in I$ with $i_0<i_1$, $l_1<l_0$ and
    $f(\a_{i_0}\b_{i_0})f(\a_{i_1}\b_{i_1}) \sim f(\a_{l_0}\b_{l_0})f(\a_{l_1}\b_{l_1})$, in
    particular
    $$f(\a_{i_0})f(\b_{i_0})f(\a_{i_1})f(\b_{i_1}) \equiv f(\a_{l_0})f(\b_{l_0})f(\a_{l_1})f(\b_{l_1})$$.


    By $\bigotimes$ \ref{equ:Phi} we have $M \models
    \ph(\a_{i_0},\b_{i_1})$ and
    $M \models \psi(\a_{l_0},\b_{l_1})$. By the assumption
    towards contradiction, $$N \models \ph^{[\frac c
    3]}(f(\a_{i_0}),f(\b_{i_1}))$$
    and $$N \models \psi^{[\frac c 3]}(f(\a_{l_0}),f(\b_{l_1}))$$
    But since
    $$f(\a_{i_0})f(\b_{i_0})f(\a_{i_1})f(\b_{i_1}) \equiv f(\a_{l_0})f(\b_{l_0})f(\a_{l_1})f(\b_{l_1})$$.

    we also get e.g.
    $$N \models \ph^{[\frac c 3]}(f(\a_{l_0}),f(\b_{l_1}))$$
    which is clearly a contradiction to the pair $(\ph,\ps)$ being
    $c$-contradictory.
\end{prf}

We conclude that relatively weak notions of categoricity imply
stability or just $c$-additive stability.



\begin{cor}\label{cor:categ}
\begin{enumerate}
\item
    Assume that there exists $\lam > |L|$ such that for every $M,N
    \in K^c$, every $\Delta \subseteq L$ finite and every $\eps>0$
    there exists $f \colon M \to N$ which is a
    $(\Delta,\eps)$-additive isomorphism. Then $\fC$ is stable.
\item
    In \emph{(i)} it is enough to assume that for every $\Delta \subseteq
    L$ finite there exists $\eps<\frac{diam(\Delta)}{3}$ and $f
    \colon M \to N$ which is a $(\Delta,\eps)$-additive isomorphism.
\item
    In \emph{(i)} it is enough to assume that for every $\Delta$ and $\eps>0$ there
    exists $\lam > |L|$ such that for every $M,N$ of density $\lam$
    there exists a $(\Delta,\eps)$-additive isomorphism from $M$ to $N$.
\item
    Similarly, it is enough to assume that for every $\Delta$ there
    exist $\eps<\frac{diam(\Delta)}{3}$ and
    $\lam > |L|$ such that for every $M,N$ of density $\lam$
    there exists a $(\Delta,\eps)$-additive isomorphism from $M$ to $N$.
\item
    Assume that for every $\Delta \subseteq
    L$ finite of diameter $\ge c$ there exists $e<\frac{c}{3}$ and $f
    \colon M \to N$ which is a $(\Delta,e)$-add-isomorphism. Then
    $\fC$ is $c$-additively stable.
\item
    Again, in (v) it is enough to assume the weaker version: for
    every $\Delta$ of diameter $\ge c$ there exists $e$ and
    $\lam > |L|$.
\end{enumerate}
\end{cor}
\begin{prf}
    Straightforward by Proposition \ref{prp:unstable} and Fact
    \ref{fct:unembed}.
\end{prf}

\section{The continuous case and Banach spaces}


The following definition of ``continuous truth'' is strongly related
to continuous model theory studied in \cite{BU} (and was, in fact, a
motivation for the second author's interest in the subject). But
since we are working here in a much more general setting than
continuous first order logic, we will not require any background or
use any results from \cite{BU} (except some notations, which we
introduce explicitly).




\begin{dfn}
\begin{enumerate}
\item
    Given a formula $\ph(\x)$ and $\a\in\fC$, we define the
    \emph{continuous truth value} of $\ph(\a)$ by
    $$\ph(\a) = \inf\set{\eps\colon \fC\models\ph^{[\eps]}(\a)}$$
    Statements of the form $[\ph(\a) = \eps]$, $[\ph(\a)\le\eps]$,
    $[\ph(\a)<\eps]$ have the obvious meaning, and we will refer
    to them as \emph{conditions}. Conditions of the form
    $[\ph(\x)=\eps]$ or $[\ph(\x)\le\eps]$ are called
    \emph{closed} while conditions of the form $[\ph(\x)<\eps]$ are
    \emph{open}.
    We will say that a tuple $\a$ \emph{satisfies} the condition
    $[\ph(\x) \le \eps]$ (or $[\ph(\x)<\eps]$, etc) if
    $[\ph(\a) \le \eps]$; sometimes we write $\a\models
    [\ph(\x)\le\eps]$.
\item
    Let $\fC$ be a momspace. We say that a formula $\ph(\x)$
    \emph{has weak negations} if for every $c>0$,
    the closed condition $[\ph(\x)\ge c]$
    is type-definable.
\end{enumerate}
\end{dfn}

\begin{rmk}
\begin{enumerate}
\item
    Note that $\ph(\a)$ is just the distance between $\a$ and the
    closed set $\ph^\fC$.
\item
    Let $M \in K$, $\ph(\x)$ a formula. Then for every
    $\a\in M$ we have $\ph^M(\a) = \ph^\fC(\a)$, that is, the
    continuous truth value does not depend on where it is being
    computed. In particular, if $N = \mcl(M) \in K^c$, then
    $\ph^M(\a) = \ph^N(\a)$ for all $\a \in M$.
\end{enumerate}
\end{rmk}

\begin{dfn}
\begin{enumerate}
\item
    Let $\ph(\x)$ be a formula, $c\le 0$. We say that $\ph$ does not have the
    \emph{$c$-additive order property} (or is \emph{$c$-additively stable},
    $\c$-add-stable)
    if for every indiscernible sequence
    $\seq{\a_i\b_i\colon i<\lam}$ long enough we have
    $$ i<j \then |\ph(\a_i,\b_j)-\ph(\a_j,\b_i)|\le c$$
\item
    Let $\ph(\x)$ be a formula, $c\ge 1$. We say that $\ph$ does not have the
    \emph{$\eps$-multiplicative order property} (or is \emph{$\eps$-multiplicatively stable},
    $\eps$-mult-stable)
    if for every indiscernible sequence
    $\seq{\a_i\b_i\colon i<\lam}$ long enough we have
    $$ i<j \then \frac 1c \le \frac{\ph(\a_i,\b_j)}{\ph(\a_j,\b_i)}\le c$$
    where we stipulate $\frac 00 = 1, \frac c0 = \infty$ for $c>0$.
\end{enumerate}
\end{dfn}

\begin{obs}\label{obs:local}
\begin{enumerate}
\item
    If a formula $\ph(\x,\y)$ is $c$-add-unstable, then $\fC$ is
    $c$-add-unstable exemplified by a pair $(\ph^{[r]}(\x,\y),\ps(\x,\y))$
    for some $r\ge 0$ and $\ps$ with $\dd(\ph,\ps)>r+c$.
\item
    $\fC$ is $c$-add-stable if and only if every formula is $c$-add-stable.
\end{enumerate}
\end{obs}
\begin{prf}
    Note that $\ph(\a,\b)>e$ means that the formula $\ph^{[e]}(\a,\b)$ does
    not hold, and therefore by the definition of a monster metric
    space there exists $\ps(\x,\y)$ such that $\ps(\a,\b)$ holds and
    the pair $(\ph^{[e]},\ps)$ is contradictory. So if e.g.
    $\ph(\a_i,\b_j)\le r$ and $\ph(\a_j,\b_i) > r+c$, then
    there exists $\ps(\x,\y)$ such that $(\ph^{[r+c]},\ps)$ is
    a contradictory pair and $\ps(\a_j,\b_i)$ holds. The rest should
    be clear.
\end{prf}

\begin{dfn}
\begin{enumerate}
\item
    Let $\ph(\x)$ be a formula, $M,N \in K$, $f \colon M\to N$, $c\ge 0$. We say that $f$ is
    a \emph{$(\ph,c)$-additive isomorphism} ($(\ph,c)$-add-isomorphism) if for
    every $\a \in M$ we have $|\ph(\a)-\ph(f(\a))|\le c$.
\item
    Let $\ph(\x)$ be a formula, $M,N \in K$, $f \colon M\to N$, $c\ge 1$. We say that $f$ is
    a \emph{$(\ph,c)$-multiplicative isomorphism} ($(\ph,c)$-mult-isomorphism) if for
    every $\a \in M$ we have $$\frac 1c \le \frac{\ph(\a)}{\ph(f(\a))}\le c$$
    where we stipulate $\frac 00 = 1, \frac c0 = \infty$ for $c>0$.
\item
    Let $\ph(\x)$ be a formula, $M,N \in K$, $f \colon M\to N$. We
    define the \emph{$\ph$-additive norm} of $f$ by
    $$\|f\|_{\ph,add} = \inf\set{\c\colon f \mbox{\;is a\;} (\ph,c)\mbox{-add-isomorhism}}$$
    (could be $\infty$).
\item
    Let $\ph(\x)$ be a formula, $M,N \in K$, $f \colon M\to N$. We
    define the \emph{$\ph$-multiplicative norm} of $f$ by
    $$\|f\|_{\ph,mult} = \inf\set{\c\colon f \mbox{\;is a\;} (\ph,c)\mbox{-mult-isomorhism}}$$
    (could be $\infty$).
\item
    Let $\ph(\x)$ be a formula, $c\ge 0$. We say that $\fC$ is
    \emph{$(\ph,c)$-additively categorical} in a cardinality $\lam$
    if for every $M,N \in K^c$ of density $\lam$ there is $f \colon
    M \to N$, $\|f\|_{\ph,add} \le c$.
\item
    Let $\ph(\x)$ be a formula, $c\ge 1$. We say that $\fC$ is
    \emph{$(\ph,c)$-multiplicatively categorical} in a cardinality $\lam$
    if for every $M,N \in K^c$ of density $\lam$ there is $f \colon
    M \to N$, $\|f\|_{\ph,mult} \le c$.
\end{enumerate}
\end{dfn}

\begin{cor}\label{cor:addcateg} Let $c\ge 0$.
\begin{enumerate}
\item
    If a formula $\ph(\x,\y)$ is $c$-add-unstable then for every
    $\lam>|L|$ there exists a sequence of models
    \lseq{M}{i}{2^\lam} of density character
    $\lam$ such that for $i\neq j$ and for every $f \colon M_i \to
    M_j$ we have $\|f\|_{\ph,add}>\frac c3$.
\item
    Assume that $\fC$ is $(\ph,c)$-additively categorical in some
    $\lam>|L|$. Then $\ph$ is $3\cdot c$-additively stable.
\item
    Assume that for every $c>0$ there is $\lam > |L|$ such that
    $\fC$ is $(\ph,c)$-additively categorical in
    $\lam$. Then $\ph$ is stable.
\end{enumerate}
\end{cor}
\begin{prf}
    This is basically a restatement of Proposition \ref{prp:unstable} and Corollary
    \ref{cor:categ}.
    Note that one has to use
    Observation \ref{obs:local} in order to obtain EM-blueprint $\Phi$ as in
    \ref{equ:Phi}.
\end{prf}

We would like to formulate the multiplicative analogue of Corollary
\ref{cor:addcateg}. It will be more convenient for us to deal with
this in the ``Hausdorff'' case. Since our main aim is connecting
weak multiplicative categoricity for classes of normed spaces with
quantifier-free stability (so $\ph(x,y) = \|x+y\|$, that is,
$\|x+y\|=0$), existence of weak negations is a reasonable
assumption.

\begin{cor}\label{cor:multcateg} Let $c\ge 1$, $\ph(\x,\y)$ a
formula with weak negations.
\begin{enumerate}
\item
    If $\ph(\x,\y)$ is $c$-mult-unstable then for every
    $\lam>|L|$ there exists a sequence of models
    \lseq{M}{i}{2^\lam} of density character
    $\lam$ such that for $i\neq j$ and for every $f \colon M_i \to
    M_j$ we have $\|f\|_{\ph,mult}>\sqrt{c}$.
\item
    Assume that $\fC$ is $(\ph,c)$-multiplicatively categorical in some
    $\lam>|L|$. Then $\ph$ is $c^2$-mult-stable.
\item
    Assume that for every $c>1$ there is $\lam > |L|$ such that
    $\fC$ is $(\ph,c)$-mult-categorical in
    $\lam$. Then $\ph$ is stable.
\end{enumerate}
\end{cor}
\begin{prf}
    (ii) and (iii) clearly follow from (i).

    (i) We just have to make sure that the proof of Proposition
    \ref{prp:unstable} can be adjusted. So suppose $\ph$ has
    the $c$-multiplicative order property, that is, there are
    arbitrary long indiscernible sequences $\seq{\a_i\colon
    i<\theta}$ satisfying e.g. for some $r$ and $\eps>0$
    $$i<j \then \ph(\a_i,\a_j) \le r, \ph(\a_j,\a_i) \ge r\cdot c+\eps$$
    (for simplicity of notation, we add dummy variables in order to combine two sequences into
    one; of course, this can be avoided, as in the proof of Proposition \ref{prp:unstable}).

    By the assumption on $\ph$, we can find an EM-blueprint such that
    the inequalities above are satisfied by the elements of the
    skeleton of $\EM(I,\Phi)$ for every linear order $I$ (so $\ps$
    in \ref{equ:Phi} is replaced with the ``weak negation''
    $[\ph(\x,\y) \ge r\cdot c+\eps]$).


    Now repeating the proof of Proposition \ref{prp:unstable}, given
    $\lam>|L|$ we get
    a sequence \lseq{M}{\al}{2^\lam} of models in $K^c$ such that

    \begin{itemize}
    \item
        $M_\al = \overline{\EM(I_\al,\Phi)}$
    \item
        $\al \neq \be \then I_\al$ is not $\theta$-embeddable into
        $I_\be$ ($\theta$ as in \ref{prp:unstable}).
    \end{itemize}

    Let $f \colon M_\al \to M_\be$. Assuming $\|f\|_{\ph,mult} \le
    \sqrt{c}$, we will get $i_0<i_1, l_1<l_0$ such that
    \begin{itemize}
    \item
        $M_\al \models [\ph(\a_{i_0},\a_{i_1}) \le r]$ and therefore
        $M_\be \models [\ph(f(\a_{i_0}),f(\a_{i_1})) \le r\sqrt{c}]$
    \item
        $M_\al \models [\ph(\a_{l_0},\a_{l_1}) \ge rc+\eps]$ and therefore
        $M_\be \models [\ph(f(\a_{l_0}),f(\a_{l_1})) \ge r\frac{c}{\sqrt{c}}+\frac{\eps}{\sqrt{c}}]$
    \item
        $f(\a_{i_0})f(\a_{i_1}) \equiv f(\a_{l_0})f(\a_{l_1})$
    \end{itemize}
    which when put together is clearly a contradiction.
\end{prf}

The following particular case is of especial interest to us in this
context.

\begin{dfn}
    Assume that $\fC$ is a \emph{normed} structure (over $\mathbb F = \setR$ or $\setC$, maybe with
    extra-structure).
    We say that $\fC$ is \emph{$c$-categorical} in a cardinal $\lam$ for $c\ge
    1$ if 
    for every two models $M$, $N$ of density
    $\lam$ there exists a linear isomorphism $f \colon M \to N$ such that for every $a
    \in M$ we have $\frac 1c \le \frac{\|a\|_M}{\|f(a)\|_N}\le c$
    (where $\frac 00 = 1, \frac r0 = \infty$ for $r \neq 0$).
\end{dfn}

\begin{rmk}
    Note that if $f$ is as in the definition above, it is a
    $(\ph(\x),c)$ - multiplicative isomorphism for any $\ph(\x) = \ph(x_{<k})$
    of the form $\ph(\x) = \|\sum_{i<k}r_ix_i\|$ (where $r_i \in
    \mathbb F$).
\end{rmk}

So by Corollary \ref{cor:multcateg} and the previous remark we
obtain the following:

\begin{cor}\label{cor:banach stable}
    Let $\fC$ be a normed structure over $\mathbb F$ (possibly with extra-structure),
    $c$-categorical in $\lam>|L|$. Then every formula
    $\ph(\x,\y) = \ph(x_{<k},y_{<\ell})$ of the form
    $$\ph(\x,\y) = \left\|\sum_{i<k} r_ix_i+\sum_{j<\ell}s_jy_j\right\|$$
    is $c^2$-mult-stable.
\end{cor}

\section{Junge-Rosenthal stability}

After hearing the statement of Corollary \ref{cor:banach stable},
Ward Henson asked whether the proof can be modified in order to
obtain a stronger version of stability. In this section we give a
positive answer to Henson's question.

The following definition in the context of Banach spaces is due to
Junge and Rosenthal (although to the best of our knowledge the
original paper was never finished). One reference is \cite{JKO}. We
will not make use of the name ``asymptotically symmetric'' suggested
in \cite{JKO} (since the definition below is simply a very natural
generalization of stability).

For the sake of simplicity of presentation, we only deal with the
multiplicative case (and assume existence of weak negations). The
additive analogue can be developed similarly.

\begin{dfn}
    Let $\ph(\x_0, \ldots, \x_{n-1})$ be a formula with weak negations,
    $c \ge 1$.
    We say that $\ph(\x_0,\ldots,\x_{n-1})$ is $(c,n)$-multiplicatively stable ($(c,n)$-mult-stable
    or just $c$-mult-stable, since $n$ is clear from $\ph$)
    if for every indiscernible sequences $\seq{\a_{i,0}\ldots\a_{i,n-1}\colon
    i<\om}$ and permutations $\sigma, \pi$ of $n$ we have
    $$ \frac 1c \le \frac{\ph(\a_{\sigma(0),0}\ldots\a_{\sigma(n-1),n-1})}{\ph(\a_{\pi(0),0}\ldots\a_{\pi(n-1),n-1})} \le c $$
\end{dfn}

The following is a modification of Fact \ref{fct:unembed}.

    \begin{fct}\label{fct:unembed_n}
    Let $\sigma$ be a permutation of $n$ and let
    $\theta_\sigma(x_0\ldots x_{n_1},y_0\ldots y_{n-1})$ be the following formula in the language of
    order: $$\theta_\sigma = (x_0<x_1<\ldots<x_{n-1})\land (y_{\sigma(0)}<y_{\sigma(1)}<\ldots<y_{\sigma(n-1)})$$
    Let
    $\lam>|L|$ regular, $\mu\ge\lam$. Then there exists a family $\lseq{I}{\al}{2^\lam}$
    of linear orders of cardinality $\mu$ which are pairwise
    $\theta_\sigma(\x,\y)$-unembeddable into each other.
\end{fct}
\begin{prf}
    By Claim 2.29 in \cite{Sh:e}, III.
\end{prf}

\begin{thm}\label{thm:JRcateg} Let $c\ge 1$, $\ph(\vec{x}) = \ph(\x_0, \ldots,\x_{n-1})$ a
formula with weak negations.
\begin{enumerate}
\item
    If $\ph(\vec{x})$ is $c$-mult-unstable then for every
    $\lam>|L|$ regular and every $\mu\ge\lam$ there exists a sequence of models
    \lseq{M}{i}{2^\lam} of density character
    $\mu$ such that for $i\neq j$ and for every $f \colon M_i \to
    M_j$ we have $\|f\|_{\ph,mult}>\sqrt{c}$.
\item
    Assume that $\fC$ is $(\ph,c)$-multiplicatively categorical in some
    $\lam>|L|$. Then $\ph$ is $c^2$-mult-stable.
\item
    Assume that for every $c>1$ there is $\lam > |L|$ such that
    $\fC$ is $(\ph(\vec{x}),c)$-mult-categorical in
    $\lam$. Then $\ph(\vec{x})$ is stable (that is, $(n,1)$-stable).
\end{enumerate}
\end{thm}
\begin{prf}
    As in the proof of Corollary \ref{cor:multcateg},
    (ii) and (iii) follow from (i), whereas for (i) we should modify the proof of
    Corollary
    \ref{cor:multcateg} appropriately. So suppose $\ph(\vec{x})$ is $c$-mult-unstable, that is, there are
    arbitrary long indiscernible sequences $\seq{\a_{i,0}\ldots\a_{i,n-1}\colon
    i<\theta}$ satisfying e.g. for some $r$, $\eps>0$ and a
    permutation $\sigma$ of $n$
    $$i_0<i_2<\ldots<i_{n-1} \then \ph(\a_{i_0,0},\ldots,\a_{i_{n-1},n-1}) \le
    r, \;\;
    \ph(\a_{i_{\sigma(0)},0},\ldots,\a_{i_{\sigma(n-1)},n-1}) \ge r\cdot c+\eps$$

    We can find an EM-blueprint such that
    the inequalities above are satisfied by the elements of the
    skeleton of $\EM(I,\Phi)$ for every linear order $I$.


    Now
    applying Fact \ref{fct:unembed_n}, given
    $\lam>|L|$ we get
    a sequence \lseq{M}{\al}{2^\lam} of models in $K^c$ such that

    \begin{itemize}
    \item
        $M_\al = \overline{\EM(I_\al,\Phi)}$
    \item
        $\al \neq \be \then I_\al$ is not $\theta_\sigma$-embeddable into
        $I_\be$
    \end{itemize}

The rest is exactly as in the proof of Corollary
\ref{cor:multcateg}.
\end{prf}

\begin{rmk}
    Note that Theorem \ref{thm:JRcateg} is not quite a generalization
    of Corollary \ref{cor:multcateg}: here we only get $2^\lam$
    models of density $\lam$ if $\lam$ is \emph{regular} (otherwise,
    we get $2^\mu$ for every regular $\mu<\lam$), whereas Corollary
    \ref{cor:multcateg} gives the maximal number in every
    $\lam>|L|$. The reason is that Theorem \ref{thm:JRcateg} relies
    on an easier straightforward argument given in \cite{Sh:e}, chapter III,
    whereas in Corollary \ref{cor:multcateg} we could apply a more
    sophisticated analysis of \cite{Sh:e}, chapter VI.
\end{rmk}

Let us conclude with a precise statement of Theorem \ref{thm:assym},
which we think of as one of the main results of the paper. It
follows immediately from Theorem \ref{thm:JRcateg}.

\begin{cor}\label{cor:banach JRstable}
    Let $\fC$ be a normed structure over $\mathbb F$ (possibly with extra-structure),
    $c$-categorical in $\lam>|L|$. Then $\fC$ is $c^2$-asymptotically
    symmetric, that is, for every $n$, the formula
    $$\ph(x_1,\ldots,x_n) = \|x_1+\ldots+x_n\|$$
    is $c^2$-mult-stable.
\end{cor}

\bibliography{common.bib}
\bibliographystyle{alpha}
\end{document}